\documentclass[journal,twoside]{IEEEtran}
\usepackage{cite}
\usepackage[cmex10]{amsmath}
\usepackage{amssymb}
\usepackage{mathtools}
\usepackage{graphicx}
\usepackage{multirow}
\usepackage{booktabs}
\usepackage{color}
\usepackage{algorithm}
\usepackage{algorithmicx,algpseudocode}
\usepackage{url}
\usepackage{float}
\usepackage{textcomp}
\allowdisplaybreaks[4]
\newcommand{\beq}{\begin{equation}}
\newcommand{\eeq}{\end{equation}}
\newcommand{\bsq}{\begin{subequations}}
\newcommand{\esq}{\end{subequations}}
\newcommand{\bq}{\begin{eqnarray}}
\newcommand{\eq}{\end{eqnarray}}
\newcommand{\bqn}{\begin{eqnarray*}}
\newcommand{\eqn}{\end{eqnarray*}}
\def\BibTeX{{\rm B\kern-.05em{\sc i\kern-.025em b}\kern-.08em
    T\kern-.1667em\lower.7ex\hbox{E}\kern-.125emX}}
\begin{document}

\title{An Improved Surrogate Method for Solving the Energy Storage Optimal Bidding Problem}
\author{Yue Chen, Wei Wei, Tongxin Li, Yunhe Hou, Feng Liu, and Jo{\~a}o P. S. Catal{\~a}o}
\maketitle

\begin{abstract}
Energy storage is expected to play an increasingly important role in mitigating variations that come along with the growing penetration of renewable energy. In this paper, we study the optimal bidding of an energy storage unit in a semi-centralized market. The energy storage unit offers its available storage capacity and maximum charging/discharging rate to the operator; then the operator clears the real-time market by minimizing the total cost. The energy storage unit is paid/charged at locational marginal price (LMP). The problem casts down to a bilevel optimization problem with a mixed-integer lower-level. An improved surrogate-based method with the combined spatial-temporal entropy term is developed to solve this problem. Numerical examples demonstrate the scalability, efficiency, and accuracy of the proposed method.
\end{abstract}


\section{Introduction}
\label{sec:introduction}
\IEEEPARstart{E}{nergy} storage (ES) can help mitigate the fluctuation of renewable energy \cite{chen2020power}. Its global capacity approximated 159 GW in 2017 and continued to grow. Fruitful researches have been conducted on the centralized management of ES units, including its optimal sizing \cite{irshad2020battery} and optimal operation \cite{abdulla2016optimal} problems. With the decentralization of the electricity market and the accelerated deployment of distributed renewable energy resources, maintaining real-time market stability by allowing ES to participate has become a crucial topic.
Currently, the regulatory framework requires that ES be treated either as a transmission asset \cite{taylor2014financial} or a wholesale market-based asset \cite{huang2017market}. 
For the latter one, there are two different market mechanisms, the \emph{semi-centralized} scheme and the \emph{deregulated} scheme. Under the semi-centralized scheme, the ES owner bids the maximum amount of energy that can be charged/ discharged and the available capacity to the operator; then the operator clears the market by solving an economic dispatch problem constrained by the received bids \cite{mohsenian2015coordinated}. Under the deregulated scheme, the ES owner has full rights to operate the storage and reference \cite{huang2017market} proved its equivalence to the semi-centralized scheme that we adopt in this paper.


This paper considers the optimal bidding problem of an energy storage unit formulated as a bilevel program with binary variables in the lower-level, making the traditional methods for solving the bilevel program not directly applicable. In fact, the proposed bilevel program can be treated as a black-box optimization problem, whose objective function value as well as the derivative information is expensive to evaluate.
Techniques to solve such a problem are heuristical methods, derivative-free methods, and surrogate-based methods. Their pros and cons are compared in \cite{vu2017surrogate}. In general, surrogate-based methods appear to have a lighter computational burden. Moreover, when applying the surrogate-based method, the ES owner only needs to know the market-clearing prices and quantities to evaluate its revenue at each sample point, which can be output by the operator without revealing other information of the market. Therefore, the surrogate-based method is more practicable than the conventional bilevel method where the ES owner requires full information about the market.

Typical surrogate-based methods are polynomials and basis function (BF) approaches \cite{vu2017surrogate}.
Polynomials cannot fit smooth functions of any shapes and can be time-consuming.
The BF approaches are better performed in many cases. This paper is based on Kriging, one of the most well-known BF approaches with a clear statistic explanation \cite{vu2017surrogate}. Vast literature focused on improving the acquisition function to enlarge the global optimization capability of the surrogate-based method \cite{regis2007stochastic, bemporad2020global}. This paper enhances the performance of the surrogate method by introducing the combined spatial-temporal (CST) Entropy, and it improves from existing work in two aspects: 1) we can directly search for the best new sample point instead of choosing from a limited set of random points. 2) the extra term in the designed acquisition function has a clear interpretation as the entropy to capture the characteristics of sample points sequence. The contributions are twofold:

1) \textbf{Improved surrogate method}. In this paper, we improve the performance of the surrogate method by proposing the CST-Entropy that can quantify the distribution of sampling points. A new acquisition function is suggested with this CST-Entropy term to facilitate the surrogate method to better search the unexplored regions. Case studies show that this method can locate an approximate optimal solution more efficiently than some renowned Kriging models and derivative-free methods, such as pattern-search and genetic algorithm (GA).

2) \textbf{Energy storage bidding model without relaxation of the constraints prohibiting simultaneous charging and discharging.} A bilevel model for energy storage bidding is established. In the upper level, the ES owner bids its maximum power and energy capacities to the operator and is paid/charged at locational marginal price (LMP). The lower-level simulates the market clearing process and outputs the LMP. The lower-level problem is a mixed-integer linear program due to the ES related constraints. Different from previous work that relaxes the constraints preventing simultaneous charging and discharging into continuous constraints \cite{li2018extended}, this paper applies the improved surrogate method to obtain the optimal strategy directly without relaxation. Therefore, our method can be applied to market clearing problems with more general objectives, where the exactness condition of the relaxation may not hold.

\section{Bilevel Energy Storage Bidding Model}
\begin{figure}[t]
\centering
\includegraphics[width=0.8\columnwidth]{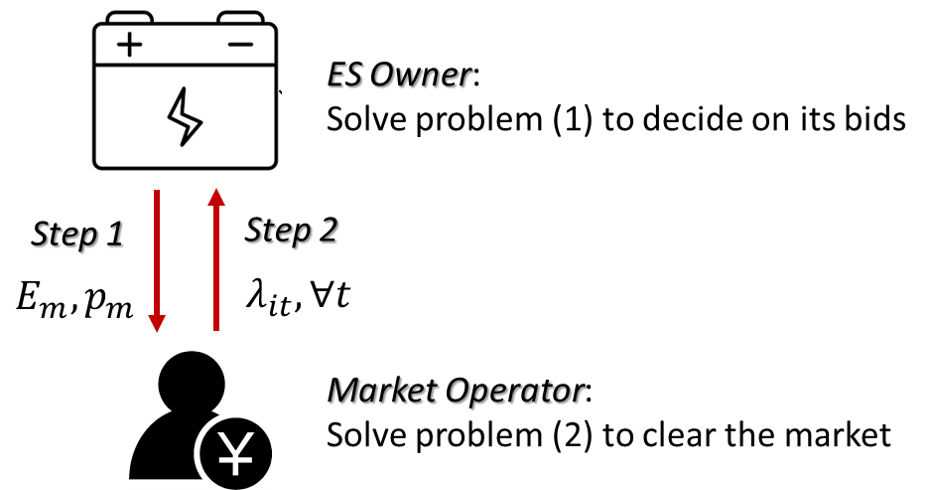}
\caption{Structure of the energy storage bidding market.}
\label{fig:market}
\end{figure}
In this section, the energy storage bidding model is given. Here, the ES runs in a semi-centralized manner \cite{huang2017market} as in Fig. \ref{fig:market}, where the state-of-charge (SoC) dynamic of ES is monitored by the market operator. First, the
ES owner offers its storage capacity $e_m$ and maximum charging$/$discharging rate $p_m$ to the market operator; then with these bids, the market operator clears the real-time market and returns the LMP to the ES owner; finally, the ES owner is charged or paid at LMP. The ES owner aims to maximize its profit \eqref{eq:ES-owner.1} by solving:
\bsq
\label{eq:ES-owner}
\begin{align}
    \mathop{\max}_{e_m,p_m} ~ & \sum \nolimits_{t=1}^T \lambda_{it: \mathbb{I}(i)=1}(p_t^d-p_t^c) \label{eq:ES-owner.1}\\
    \mbox{s.t.}~ & 0\le e_m \le E_m, ~  0 \le p_m \le P_m \label{eq:ES-owner.2}
\end{align}
\esq
For each period $t=1,...,T$, $\lambda_{it}$ is the market price given by \eqref{eq:market}; $p_t^c/p_t^d$ is the contract charging/discharging power of energy storage; $\mathbb{I}(i)=1$ if and only if the ES is connected to bus $i$. $P_m$ and $E_m$ in \eqref{eq:ES-owner.2} are the physical maximum power and energy capacities, while $p_m$ and $e_m$ are the bidding strategies. 

The market operator clears the market by solving:
\bsq
\label{eq:market}
\begin{align}
    \mathop{\min}_{p_{it},\forall i,t} ~ & \sum \nolimits_{t=1}^T \sum \nolimits_{i=1}^I \left( c_ip_{it}^2+o_ip_{it}\right) \label{eq:market.1}\\
    \mbox{s.t.}~ & y_{t+1}=y_t +\eta^c p_t^c - p_t^d / \eta^d,~\forall t \in \mathcal{T}/\{T\} \label{eq:market.2}\\
    ~ & 0 \le y_t \le e_m,~\forall t \in \mathcal{T} \label{eq:market.3}\\
    ~ & 0 \le p_t^c \le z_t^c p_m ,~ 0 \le p_t^d \le z_t^d p_m,~\forall t \in \mathcal{T} \label{eq:market.4}\\
    ~ & z_t^c, z_t^d \in \{0,1\}, z_t^c+z_t^d \le 1,~\forall t \in \mathcal{T} \label{eq:market.5}\\
    ~ & 0 \le p_{it} \le P_i,~\forall i \in \mathcal{I}, \forall t \in \mathcal{T} \label{eq:market.6}\\
    ~ & -K_m \le p_{i(t+1)}-p_{it} \le K_m,~\forall t \in \mathcal{T}/\{T\} \label{eq:market.7}\\
    ~ &  p_{it}-L_{it}=\sum \nolimits_{\mathbb{I}(i)=1} (p_t^c-p_t^d) \nonumber\\
    ~ & +\sum \nolimits_{k} b_{ik} (\theta_{it}-\theta_{kt}):\lambda_{it},\forall i \in \mathcal{I}, \forall t \in \mathcal{T}\label{eq:market.8}\\
    ~ & | b_{ik}(\theta_{it}-\theta_{kt})| \le F_{ik}  \label{eq:market.9},\forall i, k \in \mathcal{I}, \forall t \in \mathcal{T}
\end{align}
\esq
where $\mathcal{I}$ is the set of buses; $p_{it} \in [0,P_i],\forall t \in \mathcal{T}$ is the output of generator at bus $i \in \mathcal{I}$ with $c_i$, $o_i$ as its cost coefficients. For the ES, $y_t$ is the SoC; $\eta^c$/$\eta^d$, $z_t^c$/$z_t^d$ are the charging/discharging efficiency and state indicator; $K_m$ is the maximum ramping rate; $L_{it}$ refers to demand; $\theta_{it}$ is the voltage phase angle, and $b_{ik}$, $F_{ik}$ denote the susceptance and capacity of line $ik$. Constraints \eqref{eq:market.2}-\eqref{eq:market.5} are the SoC dynamics, which show the correlation between energy and power of energy storage over time. Other constraints include generator capacity \eqref{eq:market.6}, ramping limit \eqref{eq:market.7}, power balance \eqref{eq:market.8} and line flow limit \eqref{eq:market.9}. The market price is the dual variable of \eqref{eq:market.8}. The optimal bidding of ES renders a bilevel model with mixed integer linear program (MILP) in the lower level. Traditionally, bilevel program is solved by replacing the lower-level with its primal-dual optimality condition \cite{chen2018energy}, or KKT condition and linearize the complementary constraints and objective function based on the big-M method and strong duality theory \cite{fang2015coupon}. However, when the lower-level problem is an MILP, the above methods are not applicable. In the next section, we develop an improved surrogate method to solve this problem.

\section{Surrogate Method with CST-Entropy}
We develop a surrogate method with CST-Entropy to solve problem \eqref{eq:ES-owner}-\eqref{eq:market}.
Surrogate methods have been proven to be effective in solving the black-box optimization problem:
\begin{align}
\label{eq:model}
    \mathop{\min}_{x} f(x) ~ ~ \mbox{s.t.}~  l \le x \le u,   ~ x \in \mathcal{X}
\end{align}
where $f: \mathbb{R}^d \to \mathbb{R}$ is a function without an analytical form and is expensive to evaluate. $x \in \mathbb{R}^d$ is the decision variable, whose lower and upper bounds are given by $l,u \in \mathbb{R}^d$, respectively. $\mathcal{X}$ is a compact subset of $\mathbb{R}^d$, representing other constraints on $x$. Specially for problem \eqref{eq:ES-owner}-\eqref{eq:market}, $x=[e_m,p_m]^T$ with $l=[0,0]^T$ and $u=[E_m,P_m]^T$. Let $f(x)=\sum_{t=1}^T \lambda_{it:\mathbb{I}(i)=1} (p_t^d-p_t^c)$, since $\lambda_{it:\mathbb{I}(i)=1}$ and $p_t^d,p_t^c,\forall t$ are all functions of $x$, the objective $f(x)$ is also a function of $x$ solely. The value and derivative of $f(x)$ are hard to obtain as it involves solving the lower-level MILP first, so it is a black-box function.

\begin{figure}[!t]
\centering
\includegraphics[scale=0.8]{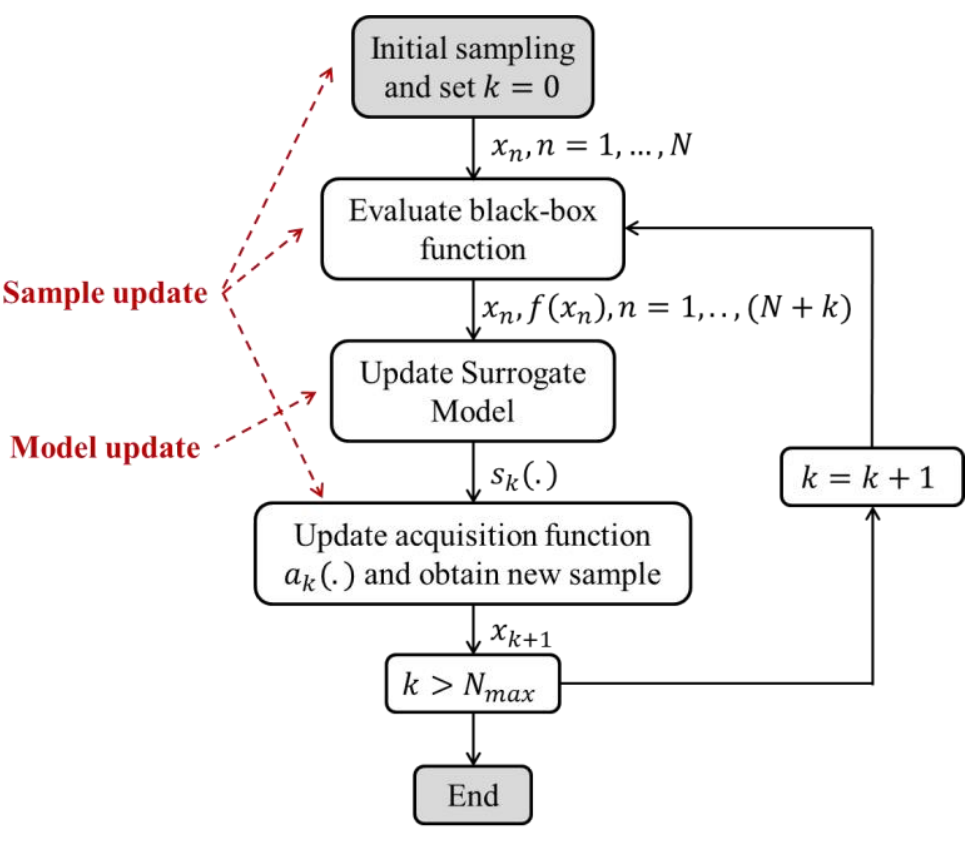}
\caption{Procedure of the surrogate optimization method.}
\label{fig:procedure}
\end{figure}

The procedure of the proposed surrogate method is shown in Fig. \ref{fig:procedure}. First, an initial set of points are given by sampling methods, e.g. Latin hypercube \cite{mckay2000comparison} used in this paper. With the samples $(x_n,f(x_n)),\forall n \in \mathcal{N}$, a surrogate model $s(x)$ is constructed to approximate $f(x)$. Then, a new sample point is found via an acquisition function $a(x)=s(x)+\xi(x)$ that trades off between \emph{exploitation} (minimizing the value of $s(x)$) and \emph{exploration} (searching unexplored regions by $\xi(x)$). Adding the new sample to the sample set, we repeat the above steps for $N_{max}$ times. Among the above procedures, the design of $a(x)$ is the key and the focus of this paper.
\subsection{Surrogate Model: Kriging}
Our model is based on Kriging. Suppose the function value $f(x)$ is a realization of a random variable $F(x)$ at point $x$, with mean $\mathbb{E}[F(x)]=\mu$ and variance $\mbox{Var}[F(x)]=\sigma^2$. Given $(x_n,f(x_n)),\forall n \in \mathcal{N}$, we can calculate the covariance matrix:
\begin{align}
    \mbox{Cov}(F)=\sigma^2R
\end{align}
Matrix $R$ is symmetric with its $(n,m)$ element given in
\begin{align}
\label{eq:corr}
    \mbox{Corr}[F(x_n),F(x_m)]=\mbox{exp} \left(-\sum\nolimits_{j=1}^d \upsilon_j |x_{nj}-x_{mj}|^{w_j}\right) \nonumber
\end{align}
Here, subscript $j$ is the index of components in vectors $x_m$ and $x_n$; $\upsilon_j$, $w_j$ are parameters. A larger $\upsilon_j$ is assigned when the function is active in the $j$-th variable. Values of $w_j$ near 2 help model smooth functions, while values of $w_j$ near 0 help model rough, non-differentiable functions. Denote $f:=\{f(x_n),\forall n\in \mathcal{N}\}$. An optimal estimation of $\mu$ and $\sigma$ via likelihood function maximization is:
\begin{equation}
\begin{aligned}
\hat \mu = \frac{1^T R^{-1} f}{1^T R^{-1} 1},~
\hat \sigma = \frac{(f-1\hat \mu )^T R^{-1} (f-1\hat \mu)}{N}
\end{aligned}
\end{equation} 
We can use parameters $\hat\mu$, $\hat \sigma$, $\upsilon_j, w_j,\forall j=1,..,d$ to predict a new point $x$. Let 
\begin{align}
    r=\left[
         \mbox{Corr}[F(x),F(x_n)], \forall n=1,...,N    \right]^T
\end{align}
be the vector of correlations between $F(x)$ and $F(x_n)$ for all $n=1,\cdots,N$.
Then the Kriging-based surrogate model is
\begin{align}
\label{eq:kriging-model}
    f(x) \approx s(x)=\hat \mu +r^T R^{-1} (f-1\hat \mu)
\end{align}

\subsection{Acquisition Function with CST-Entropy}
To obtain a new sample point, we do not simply minimize the surrogate model $s(x)$ since it may not well fit the black-box function $f(x)$. Instead, a penalty term $\xi(x)$ is included to search the unexplored region. This paper proposes a method to construct the penalty $\xi(x)$ using the idea of entropy, which can be a good measure of the dispersion of sample points. Suppose the probability of sample $x_n$ is $\pi(x_n)$, the entropy is defined as  
\begin{equation}
H(x)=-\sum \nolimits_{n=1}^N \pi(x_n) \mbox{log}(\pi(x_n))
\label{eq:Def-Entropy}
\end{equation}
However, to adapt to surrogate methods, the above formula has to be revamped: 1) Instead of probability, the location is the main concern in selecting new samples. 2) The order of sampling matters. To take these factors into account, the concept of CST-Entropy is presented, incorporating the dispersion and sequential features of sample points, denote as $\hat H(x)$.

Suppose there are $N$ existing sample points $x_n, \forall n \in \mathcal{N}$ \footnote{Since the sample point $(x_n,f(x_n))$ is uniquely determined by each $x_n$, we just write $x_n$ for short.}. When a new point $x$ is added, the incremental CST-Entropy $\Delta \hat H(x)$ can be calculated in three steps:

\textbf{Step 1: Scaling}. Both original samples $x_n,\forall n \in \mathcal{N}$ and the new point $x$ are scaled to a value between 0 and 1 through
\begin{align}
\label{eq:scale}
    \tilde{x}= (x-l)/(u-l),~ \tilde{x}_n= (x_n-l)(u-l),\forall n \in \mathcal{N}
\end{align}

\textbf{Step 2: Calculate the distance related weighting}. The Euclidean distance between $\tilde x$ and $ \tilde x_n,\forall n \in \mathcal{N}$ are given by a function $D: \mathbb{R}^d \times \mathbb{R}^d$ defined as
\begin{align}
\label{eq:distance}
    D(\tilde x,\tilde x_n)=\sqrt{(\tilde x_n-\tilde x)^T(\tilde x_n-\tilde x)},\forall n \in \mathcal{N}
\end{align}
The distance related weighting term is calculated via
\begin{align}
\label{eq:beta}
    \beta (\tilde x) =\left[\sum \limits_{n=1}^N D^{-2}(\tilde x,\tilde x_n) \right]^{-1}
\end{align}
If $\tilde x$ is close to one of the sample $\tilde x_n$, then $D^{-2}(\tilde x,\tilde x_n) \to \infty$ and therefore $\beta(\tilde x) \to 0$. In general, the function $\beta(\tilde x)$ is zero at original sample points $\tilde x_n,\forall n \in \mathcal{N}$ and grows in-between two adjacent points $\tilde x_{n_1},\tilde x_{n_2}$, $n_1, n_2 \in \mathcal{N}$.  

\textbf{Step 3: Calculate additional CST-Entropy}. The additional CST-Entropy is calculated by
\begin{align}
\label{eq:add-entropy}
    \Delta \hat H(x)=-\beta(\tilde x)\mbox{log}\beta(\tilde x)
\end{align}
where $\beta(\tilde x)$ is given by \eqref{eq:scale}-\eqref{eq:beta}. As both $\tilde x_n,\forall n \in \mathcal{N}$ and $\tilde x$ are in $[0,1]$, we have $0 \le D(\tilde x, \tilde x_n) \le 1$, so $0 \le \beta(\tilde x) \le 1$, and thus, $\Delta \hat H(x)$ is always positive.

Let $\xi(x)=-\alpha\Delta \hat H(x)$, where $\alpha$ is a parameter. The model with a larger $\alpha$ can search a broader region but takes longer to reach an optimal point, while the one with a small $\alpha$ focuses on exploitation but the search region can be limited. The acquisition function can be constructed as
\begin{align}
\label{eq:acquisition-function}
a(x)=~ \hat \mu + r^T R^{-1}(f-1\hat \mu) + \alpha\beta(\tilde x)\mbox{log}\beta(\tilde x)
\end{align}

\begin{figure}[h]
\centering
\includegraphics[width=1.0\columnwidth]{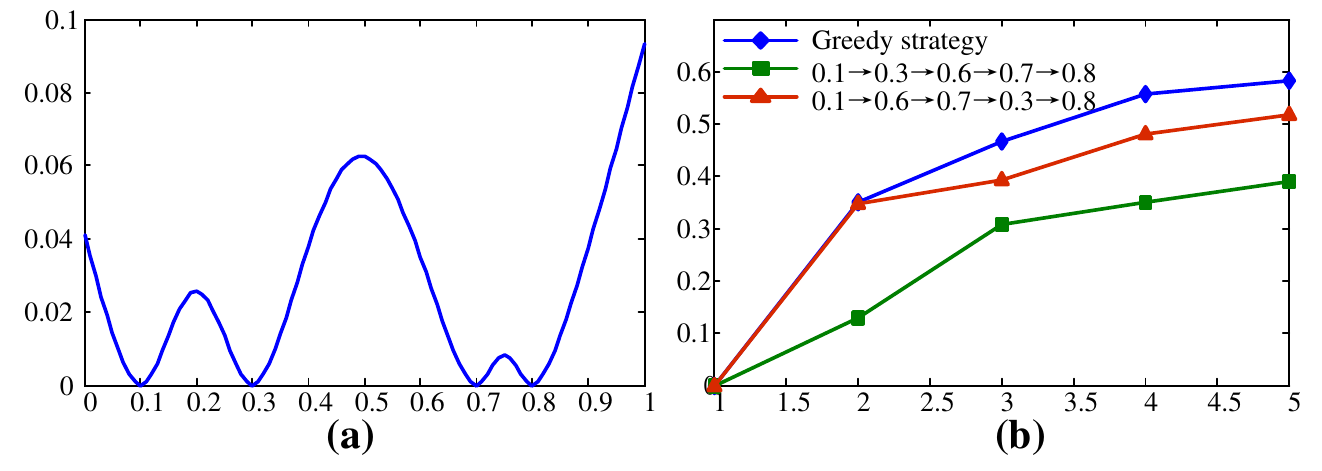}
\caption{Spatial and temporal features of CST-Entropy.}
\label{fig:entropy_geometrical}
\end{figure}

The spatial and temporal features of the CST-Entropy are illustrated in Fig. \ref{fig:entropy_geometrical}. The existing sample points are $\tilde x_1=0.1$, $\tilde x_2=0.3$, $\tilde x_3=0.7$, $\tilde x_4=0.8$. We calculate the incremental CST-Entropy \eqref{eq:add-entropy}. From Fig. \ref{fig:entropy_geometrical}(a), the values of $\Delta \hat H(x)$ are all larger than zero, meaning that the performance will not be worse-off by including one more sample. At each existing point, the $\Delta \hat H(x)$ is zero, implying that a new sample point the same as previous ones makes no improvement. The incremental CST-Entropy between two adjacent sample points becomes larger when they are farther from each other. This helps explore the regions with fewer sample points.

We then change the order following which the sampling points are added, and calculate the sum of $\Delta \hat H(x)$ step by step in Fig. \ref{fig:entropy_geometrical}(b). Different orders result in different CST-Entropy values, verifying that the order matters. It is also observed that the order generated by greedy strategy (choosing the point that maximizes \eqref{eq:add-entropy} in turns) can achieve the highest CST-Entropy value in this case. This demonstrates that if a new sample point is selected at each step by maximizing $\Delta \hat H(x)$, the sampling points will scatter more uniformly, motivating the acquisition function in the form of \eqref{eq:acquisition-function}.

\textbf{Remark:} Though a simple model \eqref{eq:ES-owner}-\eqref{eq:market} is used in this paper for better illustration, our model and method are scalable. We can consider the uncertainties by replacing problem \eqref{eq:market} with a stochastic counterpart. We can modify it to fit the day-ahead market by incorporating unit commitment constraints, where the lower-level is still an MILP and our method can still be applied. When there is more than one ES owner, we can solve each ES owner's bidding problem iteratively and obtain the equilibrium via the best-response based approach.

\section{Case studies}
\begin{figure}[t]
\centering
\includegraphics[width=1.0\columnwidth]{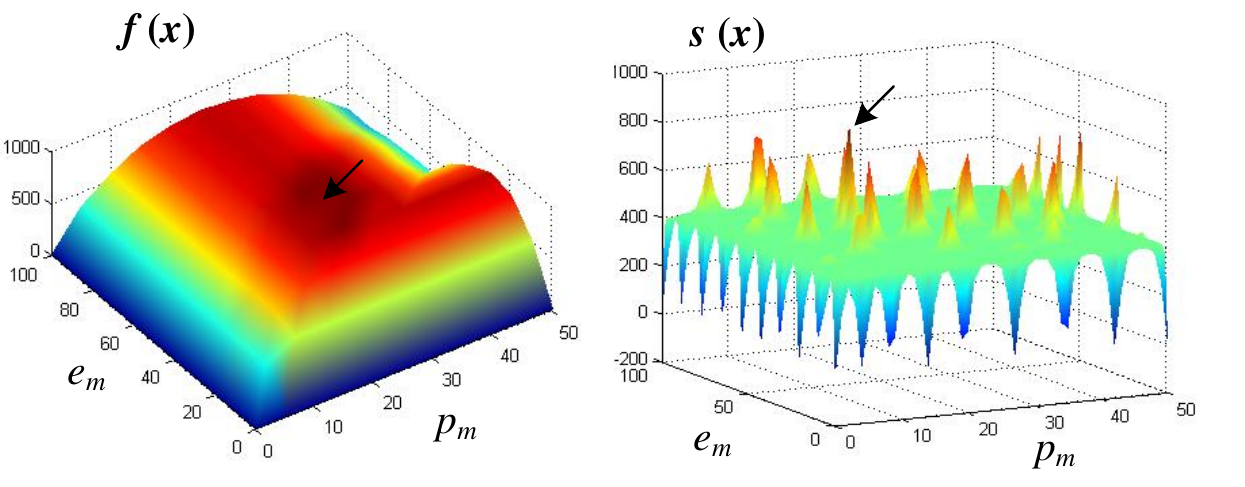}
\caption{Function $f(x)$ and Surrogate model $s(x)$ after 100 iterations.}
\label{fig:result1}
\end{figure}
\begin{figure}[t]
\centering
\includegraphics[width=0.9\columnwidth]{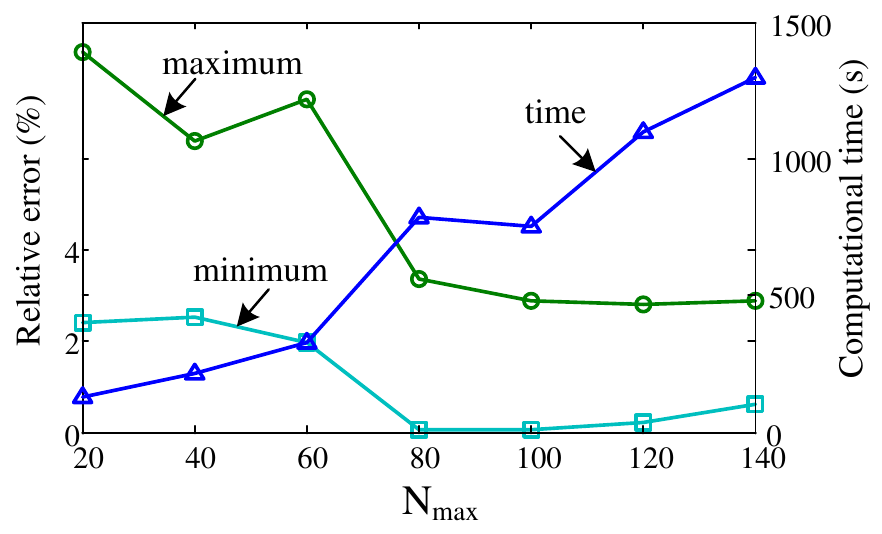}
\caption{Relative error and computational time under different $N_{max}$.}
\label{fig:result2}
\end{figure}
In this section, the performance of the proposed method is tested. We begin with a 6-bus system. Let $P_m=50$ MW, $E_m=100$ MWh, $\eta^c=\eta^d=0.8$, $T=24$ h, $N_{max}=100$; $\upsilon_j=1$, $w_j=1.5$ for all $j=1,...,d$, and $\alpha=20000$. Other data of the test systems can be found in \cite{Data}. The proposed surrogate method with CST-Entropy is applied to obtain the optimal strategy of the ES owner. The value of $s(x)$ and $f(x)$ are plotted in Fig. \ref{fig:result1}. The proposed method takes 756 seconds to find the optimal solution ($p_m=17.68$ MW, $e_m=55.84$ MWh) with a total cost of \$$885.72$. The exact optimal solution obtained by enumeration method is $p_m=17.5$ MW, $e_m=55$ MWh with a total cost of \$$886.17$. The relative error is $0.05$\%, showing the accuracy of our method.
\begin{table}[!t]
        \renewcommand{\arraystretch}{1.3}
        \centering
        \caption{Results under different $(\upsilon_j,w_j)$}
        \label{tab:parameter}
        \begin{tabular}{cccccc}
                \hline 
                Case & A & B & C & D & E  \\
                \hline
                Optimal value & 885.72 & 864.29 & 868.82 & 879.23 & 879.34\\
                Relative error & 0.05\% & 2.47\% & 1.96\% & 0.78\% & 0.77\%\\
                Time (s) & 756 & 835 & 776 & 1074 & 1847\\
                \hline
        \end{tabular}
\end{table}

We further compare the results under different $(\upsilon_j,w_j)$ as in TABLE \ref{tab:parameter}. Let Cases A-E denote the scenarios with $(\upsilon_j,w_j )=(1,1.5),(1,2),(1,0.5),(10,1.5),(0.1,1.5)$, respectively. Our method can achieve a high accuracy with all relative errors less than 3\%. The one with $(\upsilon_j,w_j)=(1,1.5)$ is the most accurate, which is also the benchmark setting in this paper. The computational times are less than 2000s, which is acceptable. We also test the impact of $N_{max}$ by changing its value from 20 to 140. We run our algorithm for 5 times with each given $N_{max}$ and the minimum/maximum relative errors are recorded in Fig. \ref{fig:result2}. The maximum relative errors are less than 8\%, and when $N_{max}$ increases, the relative errors become stable which are lower than 2.9\%. This shows that though the initial point may influence the performance of the proposed method, it is still precise enough. Moreover, the computational time increases little when $N_{max}$ grows.

Our method is also compared with some renowned derivative-free optimization methods including pattern-search and genetic algorithm (GA) in TABLE \ref{tab:comparison}, and the Kriging model with metrics response  surface  weighted  score  (MRS) in TABLE \ref{tab:Kriging}. Let $N_{max}=100$ with a time limits as $7000$ seconds. Results show that the proposed method can greatly reduce the computational time without sacrificing optimality.
\begin{table}[!t]
        \renewcommand{\arraystretch}{1.3}
        \centering
        \caption{Relative error and computational time of different methods}
        \label{tab:comparison}
        \begin{tabular}{ccccc}
                \hline 
                Method & Proposed method & Pattern-search & GA  \\
                \hline
                Relative error & 0.05\% & 0.02\% & 73.96\% \\
                Time (s) & 756 & 5929 & 7000+ \\ 
                \hline
        \end{tabular}
\end{table}
\begin{table}[!t]
        \renewcommand{\arraystretch}{1.3}
        \centering
        \caption{Comparison of different Kriging models}
        \label{tab:Kriging}
        \begin{tabular}{ccccc}
                \hline 
                 & 6-bus  & 6-bus  & 69-bus  & 69-bus \\
                 & + our method & + MRS & + our method & + MRS\\
                \hline
                Relative error & 0.05\% & 2.29\% & 0.001\% & 0.96\%\\
                 Time (s) & 756 & 403 & 998 & 383\\
                \hline
        \end{tabular}
\end{table}

\section{Conclusions}

This paper proposes an improved surrogate method to solve the optimal energy storage bidding problem, which casts as a bilevel program with a mixed-integer linear lower-level. To better explore the unsearched regions to enhance accuracy, the CST-entropy is proposed and added to the acquisition function. Compared with existing methods, the proposed surrogate method can reduce computational time while achieving high accuracy. Future research direction include acceleration methods to further speed up the algorithm and more delicate design of the CST-entropy to improve accuracy.

\end{document}